\documentclass[a4paper,reqno,12pt,final]{amsart}
\usepackage{amsmath,amsfonts,amssymb,amsthm,amstext,eucal,bbold}
\usepackage{array}
\usepackage[T1]{fontenc}
\usepackage[latin1]{inputenc}
\DeclareMathOperator{\AdS}{AdS}

\DeclareMathOperator{\dS}{dS}

\DeclareMathOperator{\End}{End}
\DeclareMathOperator{\Mat}{Mat}
\DeclareMathOperator{\Ric}{Ric}

\DeclareMathOperator{\dvol}{dvol}

\renewcommand{\d}{\partial}

\newcommand{\Cl}{\mathrm{C}\ell}
\newcommand{\GL}{\mathrm{GL}}
\newcommand{\RR}{\mathbb{R}}

\newcommand{\SO}{\mathrm{SO}}

\newcommand{\Spin}{\mathrm{Spin}}

\newcommand{\eL}{\mathcal{L}}
\newcommand{\eM}{\mathcal{M}}

\newcommand{\fS}{\mathfrak{S}}
\newcommand{\fg}{\mathfrak{g}}

\newcommand{\fp}{\mathfrak{p}}
\newcommand{\fk}{\mathfrak{k}}
\newcommand{\half}{\tfrac{1}{2}}

\newcommand{\spb}{\$}
\newcommand{\eH}{\mathcal{H}}
\newcommand{\eV}{\mathcal{V}}

\newtheorem{dfn}{Definition}

\newtheorem{prop}{Proposition}
\newtheorem{thm}{Theorem}

\newtheorem{miracle}{Minor Miracle}
\begin{document}

\title[Maximal supersymmetry in $10$ and $11$ dimensions]{Maximal
supersymmetry in ten and eleven dimensions}
\author[José Figueroa-O'Farrill]{José Figueroa-O'Farrill}
\address{Department of Mathematics and Statistics, University of
Edinburgh}
\email{j.m.figueroa@ed.ac.uk}
\begin{abstract}
  This is the written version of a talk given in Bonn on September
  11th, 2001 during a workshop on \emph{Special structures in string
  theory}.  We report on joint work in progress with George
  Papadopoulos aimed at classifying the maximally supersymmetric
  solutions of the ten- and eleven-dimensional supergravity theories
  with $32$ supercharges.
\end{abstract}
\maketitle

\section{Eleven-dimensional supergravity}

Eleven-dimensional supergravity was predicted by Nahm \cite{Nahm} and
constructed soon thereafter by Cremmer, Julia and Scherk \cite{CJS}.
We will only be concerned with the bosonic equations of motion.  The
geometrical data consists of $(M^{11}, g, F)$ where $(M,g)$ is an
eleven-dimensional lorentzian manifold with a spin structure and $F
\in \Omega^4(M)$ is a closed $4$-form.  The equations of motion
generalise the Einstein--Maxwell equations in four dimensions.  The
Einstein equation relates the Ricci curvature to the energy momentum
tensor of $F$.  More precisely, the equation is
\begin{equation}
  \label{eq:einstein}
    \Ric(g) = T(g,F)
\end{equation}
where the symmetric tensor
\begin{equation*}
  T(X,Y) = \half \left< \imath_X F, \imath_Y F\right> - \tfrac16
  g(X,Y) |F|^2~,
\end{equation*}
is related to the energy-momentum tensor of the (generalised) Maxwell
field $F$.  In the above formula, $\left<-,-\right>$ denotes the
scalar product on forms, which depends on $g$, and $|F|^2 =
\left<F,F\right>$ is the associated (indefinite) norm.  The
generalised Maxwell equations are now nonlinear:
\begin{equation}
  \label{eq:maxwell}
   d \star F = \half F \wedge F~.
\end{equation}

\begin{dfn}
  A triple $(M,g,F)$ satisfying the equations \eqref{eq:einstein} and
  \eqref{eq:maxwell} is called a \textbf{(classical) solution} of
  eleven-dimensional supergravity.
\end{dfn}

Let $\spb$ denote the bundle of spinors\footnote{As David Calderbank
  likes to remind me, there is big money in spin geometry.} on $M$.
It is a real vector bundle of rank $32$ with a spin-invariant
symplectic form $\left(-,-\right)$.  A differential form on $M$ gives
rise to an endomorphism of the spinor bundle via the composition
\begin{equation*}
  c: \Lambda T^*M  \xrightarrow{\cong} \Cl(T^*M) \to \End\spb~,
\end{equation*}
where the first map is the bundle isomorphism induced by the vector
space isomorphism between the exterior and Clifford algebras, and the
second map is induced from the action of the Clifford algebra
$\Cl(1,10)$ on the spinor representation $S$ of $\Spin(1,10)$.  In
signature $(1,10)$ one has the algebra isomorphism
\begin{equation*}
  \Cl(1,10) \cong \Mat(32,\RR) \oplus \Mat(32,\RR)~,
\end{equation*}
hence the map $\Cl(1,10) \to \End S$ has kernel.  In other words, the
map $c$ defined above involves a choice.  This comes down to choosing
whether the (normalised) volume element in $\Cl(1,10)$ acts as $\pm$
the identity.  We will assume that a choice has been made once and for
all.

\begin{dfn}
  We say that a classical solution $(M,g,F)$ is
  \textbf{supersymmetric} if there exists a nonzero spinor
  $\varepsilon \in \Gamma(\spb)$ which is parallel with respect to the
  \textbf{supercovariant connection}
  \begin{equation*}
    D : \Gamma(\spb) \to \Gamma( T^*M \otimes \spb)
  \end{equation*}
  defined, for all vector fields $X$, by
  \begin{equation*}
    D_X \varepsilon  = \nabla_X \varepsilon - \Omega_X(F)
    \varepsilon~,
  \end{equation*}
  where $\nabla$ is the spin connection and $\Omega(F):TM \to End\spb$
  is defined by
  \begin{equation*}
    \Omega_X(F) = \tfrac1{12} c(X^\flat \wedge F) - \tfrac16
    c(\imath_X F)~,
  \end{equation*}
  with $X^\flat$ the one-form dual to $X$.
\end{dfn}

A nonzero spinor $\varepsilon$ which is parallel with respect to $D$
is called a \textbf{Killing spinor}.  This is a generalisation of the
usual geometrical notion of Killing spinor (see, for example,
\cite{BFGK}).  The name is apt because Killing spinors are ``square
roots'' of Killing vectors.  Indeed, one has the following

\begin{prop}
  Let $\varepsilon_i$, $i=1,2$ be Killing spinors: $D \varepsilon_i =
  0$.  Then the vector field $V$ defined, for all vector fields $X$,
  by
  \begin{equation*}
    g(V,X) = \left( \varepsilon_1, X \cdot \varepsilon_2\right)
  \end{equation*}
  is a Killing vector.
\end{prop}

There is a vast literature on supersymmetric solutions of
eleven-dimensional supergravity, but so far very few results of a
general nature.  This problem comes down to studying the
supercovariant connection $D$.  Alas, $D$ is not induced from a
connection on the tangent bundle and in fact, it does not even
preserve the symplectic structure.  In fact, one has the following

\begin{prop}
  The holonomy of $D$ is generically $\GL(32,\RR)$.
\end{prop}

\section{Kaluza--Klein reduction and type IIA supergravity}

Suppose that $(M,g,F)$ is a classical solution of eleven-dimensional
supergravity admitting a free circle (or $\RR$) action leaving $g$ and
$F$ invariant.  Let $\xi$ denote the Killing vector generating this
action.  We will assume that $\xi$ is spacelike, so that its norm is
everywhere positive.  Let $N$ denote the space of orbits.  For
definiteness we can consider the case of a free circle action.  Let
$\pi : M \to N$ be the canonical projection sending a point in $M$ to
the (unique) orbit it belongs to.  For every $m\in M$, the tangent
space to $M$ at $m$ splits into vertical and horizontal subspaces:
\begin{equation*}
  T_m M = \eV_m \oplus \eH_m~,
\end{equation*}
where $\eV_m$ is the one-dimensional subspace spanned by $\xi(m)$ and
$\eH_m = \eV_m^\perp$ is its perpendicular complement.  The projection
$\pi_*$ defines an isomorphism $\eH_m \cong T_{\pi m} N$ and there is
a unique metric $h$ on $N$ for which this is also an isometry.

The horizontal distribution $\eH$ defines a one-form $\omega$ such
that $\ker \omega= \eH$ and normalised so that $\omega(\xi) = 1$.
Introducing a coordinate $\theta$ adapted to the circle action, we
have $\xi = \d_\theta$ and $\omega = d\theta + A$, where $A$ is a
horizontal one-form called the \textbf{RR one-form potential}.  Its
field-strength pulls back to the curvature $d\omega$ of the principal
connection, which is both horizontal and invariant, hence basic.

Finally, the metric on the fibres is described by a function $\Phi$ on
$N$, called the \textbf{dilaton}.  In terms of these data, the
eleven-dimensional metric can be written as
\begin{equation}
  \label{eq:KKmetric}
  g = \pi^* h + e^{\pi^* \Phi} \omega \otimes \omega~.
\end{equation}

Similarly we can decompose the four-form $F$ as follows
\begin{equation*}
  F = \omega \wedge \imath_\xi F + K~,
\end{equation*}
where $\imath_\xi K = 0$.  It follows from the fact that $F$ is closed
and invariant, that $\imath_\xi F$ and $K$ are basic.  Therefore there
are forms $H\in\Omega^3(N)$ and $G\in\Omega^4(N)$ on $N$ such that
\begin{equation*}
  \imath_\xi F = \pi^* H \qquad\text{and}\qquad K = \pi^* G~,
\end{equation*}
whence
\begin{equation}
  \label{eq:KKfourform}
  F = \omega \wedge \pi^* H + \pi^* G~.
\end{equation}
The closed $3$-form $H$ is called the \textbf{NSNS $3$-form} and $G$
is called the \textbf{RR $4$-form field-strength}.

The data $(N,h,\Phi,H,A,G)$ is then a solution of the equations of
motion of ten-dimensional type IIA supergravity theory. These
equations are obtained from equations \eqref{eq:einstein} and
\eqref{eq:maxwell} by simply inserting the expressions
\eqref{eq:KKmetric} for the metric and \eqref{eq:KKfourform} for the
four-form.

The data $(N,h,\Phi,H)$ defines the \textbf{common sector} of type II
supergravity in ten dimensions.  In this context, the NSNS $3$-form
$H$ can be interpreted as the torsion three-form of a metric
connection on $(N,h)$.  This gives rise to a variety of torsioned
geometries discussed at this conference by Friedrich and Papadopoulos.
 
How about supersymmetry?  The circle action lifts to an action on the
spinor bundle, which is infinitesimally generated by the
\textbf{spinorial Lie derivative} introduced by Lichnerowicz.  If
$\varepsilon$ is any spinor, then
\begin{equation*}
  \eL_\xi \varepsilon = \nabla_\xi \varepsilon + \tfrac14
  c(d\xi^\flat) \varepsilon~.
\end{equation*}
An invariant Killing spinor
\begin{equation*}
  D_X \varepsilon = 0 \qquad \text{and}\qquad \eL_\xi \varepsilon = 0
\end{equation*}
gives rise to a IIA Killing spinor, and viceversa (at least locally).
Notice that the IIA Killing spinor equation has a purely algebraic
component, namely
\begin{equation*}
  \left( \eL_\xi - \nabla_\xi \right) \varepsilon = 0~,
\end{equation*}
called the \textbf{dilatino equation}.

A useful principle in this game is the fact that supersymmetric
solutions to IIA supergravity can be lifted to invariant
supersymmetric solutions of eleven-dimensional supergravity.  This
procedure does not involve any loss of supersymmetry; although it may
sometimes result in accidental supersymmetry in eleven dimensions.
This means that it is often more convenient to work with
eleven-dimensional supergravity than with IIA supergravity.

\section{Maximal supersymmetry}

\begin{dfn}
  A classical solution of eleven-dimensional or type IIA supergravity
  is called \textbf{maximally supersymmetric} if the space of Killing
  spinors is of maximal dimension, namely $32$.
\end{dfn}

If $(M,g,F)$ is a maximally supersymmetric classical solution of
eleven-dimensional supergravity, the supercovariant connection $D$ is
flat.  Solving the flatness equations of the supercovariant connection
one arrives at the following theorem.

\begin{thm}[\cite{KG,FOPMax}]
  \label{thm:class}
  Let $(M,g,F)$ be a maximally supersymmetric solution of
  eleven-dimensional supergravity.  Then $(M,g)$ has constant scalar
  curvature $s$, and depending on the value of $s$ one has the
  following classification:
  \begin{itemize}
  \item If $s>0$, then $(M,g)$ is locally isometric to $\AdS_7 \times
    S^4$, where $\AdS_7$ is the lorentzian space-form of constant
    negative curvature $-7s$ and $S^4$ is the round sphere with
    constant positive curvature $8s$; and $F = \sqrt{6 s} \dvol(S^4)$.
  \item If $s<0$, then $(M,g)$ is locally isometric to $\AdS_4 \times
    S^7$, where $\AdS_4$ has constant negative curvature $8s$ and
    $S^7$ is the round sphere with constant positive curvature $-7s$;
    and $F = \sqrt{-6 s} \dvol(\AdS_4)$.
  \item If $s=0$ there are two possibilities:
    \begin{itemize}
    \item $(M,g)$ is flat and $F=0$; or
    \item $(M,g)$ is locally isometric an indecomposable lorentzian
      symmetric space with solvable transvection group, and $F\neq 0$.
    \end{itemize}
  \end{itemize}
\end{thm}

The classification of symmetric spaces in indefinite signature is
hindered by the fact that there is no splitting theorem saying that if
the holonomy representation is reducible, the space is locally
isometric to a product.  In fact, local splitting implies both
reducibility \emph{and} a nondegeneracy condition on the factors
\cite{Wu}.  This means that one has to take into account reducible yet
indecomposable holonomy representations.  The general semi-riemannian
case is still open, but indecomposable lorentzian symmetric spaces
were classified by Cahen and Wallach \cite{CahenWallach} more than
thirty years ago.  At least for dimension $n\geq 3$, there are three
types of indecomposable lorentzian symmetric spaces:
\begin{itemize}
\item $\dS_n$ (de~Sitter space), the space form with constant positive
  curvature,
\item $\AdS_n$ (anti~de~Sitter space), the space form with constant
  negative curvature, and
\item an ($n-3$)-dimensional family of ``pp-waves'' with solvable
  transvection group.
\end{itemize}
It is precisely this last class of symmetric spaces which, for $n=11$,
describes the gravitational part of a maximally supersymmetric
solution of eleven-dimensional supergravity.

\section{The Cahen--Wallach pp-waves}

The Cahen--Wallach $n$-dimensional pp-waves are constructed as
follows.  Let $V$ be a real vector space of dimension $n-2$ endowed
with a euclidean structure $\left<-,-\right>$.  Let $V^*$ denote its
dual. Let $Z$ be a real one-dimensional vector space and $Z^*$ its
dual.  We will identify $Z$ and $Z^*$ with $\RR$ via canonical dual
bases $\{e_+\}$ and $\{e_-\}$, respectively.  Let $A \in S^2V^*$ be a
symmetric bilinear form on $V$.  Using the euclidean structure on $V$
we can associate with $A$ an endomorphism of $V$ also denoted $A$:
\begin{equation*}
 \left< A(v), w\right> = A(v,w)\qquad\text{for all $v,w\in V$.}
\end{equation*}
We will also let $\flat:V \to V^*$ and $\sharp: V^* \to V$ denote the
musical isomorphisms associated to the euclidean structure on $V$.

Let $\fg_A$ be the Lie algebra with underlying vector space $V \oplus
V^* \oplus Z \oplus Z^*$ and with Lie brackets
\begin{equation}
  \label{eq:gA}
  \begin{aligned}[m]
    [e_-, v] &= v^\flat\\
    [e_-, \alpha] &= A(\alpha^\sharp)\\
    [\alpha, v] &= A(v,\alpha^\sharp) e_+~,
  \end{aligned}
\end{equation}
for all $v\in V$ and $\alpha \in V^*$.  All other brackets not
following from these are zero.  The Jacobi identity is satisfied by
virtue of $A$ being symmetric.  Notice that since its second derived
ideal is central, $\fg_A$ is (three-step) solvable.

Notice that $\fk_A=V^*$ is an abelian Lie subalgebra, and its
complementary subspace $\fp_A= V \oplus Z \oplus Z^*$ is acted on by
$\fk_A$.  Indeed, it follows easily from \eqref{eq:gA} that
\begin{equation*}
  [\fk_A,\fp_A] \subset \fp_A \qquad\text{and}\qquad [\fp_A,\fp_A]
  \subset \fk_A~,
\end{equation*}
whence $\fg_A = \fk_A \oplus \fp_A$ is a symmetric split.  Lastly, let
$B\in\left(S^2\fp_A^*\right)^{\fk_A}$ denote the invariant symmetric
bilinear form on $\fp_A$ defined by
\begin{equation*}
  B(v,w) = \left<v,w\right> \qquad\text{and}\qquad B(e_+,e_-)=1~,
\end{equation*}
for all $v,w\in V$.  This defines on $\fp_A$ a $\fk_A$-invariant
lorentzian inner product of signature $(1,n-1)$.

We now have the required ingredients to construct a (lorentzian)
symmetric space.  Let $G_A$ denote the connected, simply-connected Lie
group with Lie algebra $\fg_A$ and let $K_A$ denote the Lie subgroup
corresponding to the subalgebra $\fk_A$.  The lorentzian inner product
$B$ on $\fp_A$ induces a lorentzian metric $g$ on the space of cosets
\begin{equation*}
  M_A = G_A/K_A~,
\end{equation*}
turning it into a symmetric space.

\begin{prop}[\cite{CahenWallach}]
  The metric on $M_A$ defined above is indecomposable if and only if
  $A$ is nondegenerate.  Moreover, $M_A$ and $M_{A'}$ are isometric if
  and only if $A$ and $A'$ are related in the following way:
  \begin{equation*}
    A'(v,w) = c A(Ov, Ow) \qquad\text{for all $v,w\in V$,}
  \end{equation*}
  for some orthogonal transformation $O:V\to V$ and a positive scale
  $c>0$.
\end{prop}

From this result one sees that the moduli space $\eM_n$ of
indecomposable such metrics in $n$ dimensions is given by
\begin{equation*}
  \eM_n = \left( S^{n-3} - \Delta \right) / \fS_{n-2}~,
\end{equation*}
where
\begin{equation*}
  \Delta = \left\{ (\lambda_1,\dots,\lambda_{n-2}) \in S^{n-3} \subset
    \RR^{n-2} \mid \lambda_1 \cdots \lambda_{n-2} = 0\right\}
\end{equation*}
is the singular locus consisting of eigenvalues of degenerate
$A$'s, and $\fS_{n-2}$ is the symmetric group in $n-2$ symbols, acting
by permutations on $S^{n-3} \subset \RR^{n-2}$.

A remarkable fact which is still not properly understood is the
following

\begin{miracle}
  There is a unique point $A_* \in \eM_{11}$ for which $(M_{A_*}, g)$
  is the gravitational part of a maximally supersymmetric solution of
  eleven-dimensional supergravity.
\end{miracle}

Explicitly, we can write this solution as
\begin{equation*}
  \begin{aligned}[m]
    g &= 2 dx^+ dx^- - \left(\sum_{i=1}^3 (x^i)^2 + \tfrac14
      \sum_{i=4}^9  (x^i)^2 \right) \left(dx^-\right)^2 + \sum_{i=1}^9
    \left(dx^i\right)^2\\
    F &= 3 dx^- \wedge dx^1 \wedge dx^2 \wedge dx^3~.
  \end{aligned}
\end{equation*}

The isometry group of the metric is not just $G_A$ but the larger
group
\begin{equation*}
  G_A \rtimes \left( \SO(3) \times \SO(6) \right)~,
\end{equation*}
where $\SO(3) \times \SO(6) \subset \SO(9)$ acts on $G_A$ by
exponentiating the restriction of the natural action of $\SO(9)$ on $V
\oplus V^* = \RR^9 \oplus \RR^9$.  Intriguingly, the dimension of the
isometry group is $38$, which is the same as the dimension of the
isometry groups of the other maximally supersymmetric solutions of
AdS-type.  This deserves to be better understood.

\section{Maximal supersymmetry in type IIA supergravity}

Let $(N,h,\Phi,H,A,G)$ be a maximally supersymmetric solution of type
IIA supergravity.  Let $(M,g,F)$, where $g$ and $F$ are given by
\eqref{eq:KKmetric} and \eqref{eq:KKfourform} respectively, denote the
corresponding circle-invariant\footnote{This is a local result -- the
  action is only infinitesimal, hence we cannot distinguish between
  circle or $\RR$ actions.} solution of eleven-dimensional
supergravity.  Since no supersymmetry is lost in this process,
$(M,g,F)$ is also maximally supersymmetric.  Moreover, the action of
the Killing vector $\xi$ must leave invariant \emph{all} Killing
spinors.

It is then a matter of going through the maximally supersymmetric
solutions classified in Theorem~\ref{thm:class} and checking whether
there exists a Killing vector which leaves all Killing spinors
invariant.  For the AdS solutions, it follows from the semisimplicity
of the isometry algebra that no such Killing vector exists.  It was
shown in \cite{FOPflux}, albeit in a different context, that neither
does the maximally supersymmetric pp-wave solution admit such Killing
vectors.  Finally, for the flat solution with $F=0$, we can let $\xi$
be any translation along a spacelike direction.  The resulting IIA
solution is such that $(N,h)$ is flat, the dilaton is constant and all
other fields vanish.  In summary, we have proven the following.

\begin{thm}
  The only maximally supersymmetric solution of type IIA supergravity
  is a flat spacetime with constant dilaton and vanishing $(A,H,G)$.
\end{thm}

Since only the common sector fields are nonzero, this solution is also
a maximally supersymmetric solution of type IIB supergravity.  However
in this case we know at least another class of maximally
supersymmetric solutions, with geometry $\AdS_5 \times S^5$.  The
classification of maximally supersymmetric solutions of type IIB
supergravity is work in progress \cite{FOPMax}.

\section*{Acknowledgments}

It is a pleasure to thank George Papadopoulos for the ongoing
collaboration on this project, and David Calderbank and Michael Singer
for conversations.  I would like to express my thanks to Dmitri
Alekseevsky, Vicente Cortés, Chand Devchand and Toine Van Proeyen for
the invitation to participate in the workshop, and to the DFG for
their financial support.  I am a member of EDGE, Research Training
Network HPRN-CT-2000-00101, supported by The European Human Potential
Programme.

%

\providecommand{\bysame}{\leavevmode\hbox to3em{\hrulefill}\thinspace}

\end{document}